\documentclass[a4paper,reqno,10pt]{amsart}

\usepackage{geometry,times}
\usepackage{amsmath,amssymb}

\newtheorem{Theorem}{Theorem}[section]

\newtheorem{Lemma}[Theorem]{Lemma}

\newtheorem{Algorithm}[Theorem]{Algorithm}
\newenvironment{Proof}
{\begin{trivlist}\item[]{{\sc Proof.}}}{\hfill{$\square$}\noindent\end{trivlist}}

\begin{document}
\title[On the generation of Heronian triangles]{On the generation of Heronian triangles}

\author{Sascha Kurz}

\address{Sascha Kurz\\Department of Mathematics, Physic and Informatics\\University of Bayreuth\\Germany}
\email{sascha.kurz@uni-bayreuth.de}

\begin{abstract}
  We describe several algorithms for the generation of integer Heronian triangles
  with diameter at most $n$. Two of them have running time $\mathcal{O}\left(n^{2+\varepsilon}\right)$.
  We enumerate all integer Heronian triangles for $n\le 600000$ and apply the complete
  list on some related problems.
\end{abstract}


\maketitle

\section{Introduction}

  \noindent
  The greek mathematician Heron of Alexandria (c. 10 A.D. -- c. 75 A.D.) was probably the first to prove a
  relation between the side lengths $a$, $b$, and $c$  and the area $A$ of a triangle,
  $$
    A=\sqrt{s(s-a)(s-b)(s-c)}\quad\quad\mbox{where }s=\frac{a+b+c}{2}\,.
  $$
  If its area and its side lengths are rational then it is called a Heronian triangle.\\
  Triangles with integer sides and rational area were considered by the Indian mathematician Brahmagupta 
  (598-668 A.D.) who gives the parametric solution 
  \begin{eqnarray}
    a&=&\frac{p}{q}h(i^2+j^2)\nonumber\\
   b&=&\frac{p}{q}i(h^2+j^2)\nonumber\\
   c&=&\frac{p}{q}(i+h)(ih-j^2)\nonumber
  \end{eqnarray} 
  for positive integers $p$, $q$, $h$, $i$, and $j$ satisfying $ih>j^2$ and $\gcd(p,q)=\gcd(h,i,j)=1$.

  Much has been contributed \cite{0923.11052,0219.10027,0893.11009,0986.51030,0255.10019,0937.11012} 
  to the determination of such integral triangles, but still little is known about
  the generation of integer Heronian triangles with diameter $n=\max(a,b,c)$.
  Our aim is to develop a fast algorithm for the generation of the complete set of integer Heronian
  triangles with diameter $n$.

  In this context an extensive search on those triangles was made by Randall L. Rathbun \cite{petersen}. He simply checked the
  $7818928282738$ integer triangles with diameter at most $2^{17}$ and received
  $5801746$ primitive, i.~e. such with $\gcd(a,b,c)=1$, integer Heronian triangles with rational area.

  In the next section we introduce a new parameterization and in Section \ref{sec_algo} 
  we give some algorithms for the generation of integer Heronian triangles. We finish with some combinatorial problems 
  connected to Heronian triangles.

\section{A new parameterization}

\noindent
The obstacle for a computational use of Brahmagupta's parametric solution is the denominator $q$. So we first prove
a few lemmas on $q$.

\begin{Lemma}
  We can assume that the denominator $q$ can be written as $q=w_1w_2w_3w_4$ with pairwise coprime integers $w_1,w_2,w_3,w_4$ and
 \begin{tabbing}
    \hspace{25mm}\=\hspace{25mm}\=\hspace{36mm}\=\kill
    $(w_1,h)=w_1$, \> $(w_1,i)=w_1$,  \> $(w_1,i^2+j^2)=1$, \> $(w_1,h^2+j^2)=1$,   \\
    $(w_2,h)=w_2$, \> $(w_2,i)=1$,    \> $(w_2,i^2+j^2)=1$, \> $(w_2,h^2+j^2)=w_2$, \\
    $(w_3,h)=1$,   \> $(w_3,i)=w_3$,  \> $(w_3,i^2+j^2)=w_3$, \> $(w_3,h^2+j^2)=1$, \\
    $(w_4,h)=1$,   \> $(w_4,i)=1$,    \> $(w_4,i^2+j^2)=w_4$, \> $(w_4,h^2+j^2)=w_4$,
  \end{tabbing} 
  where $(x,y)$ abbreviates $\gcd(x,y)$.
\end{Lemma}
\begin{Proof}
  Suppose $q=\frac{q_1q_2}{\gcd(q_1,q_2)}$ with $q_1|h$ and $q_2|i^2+j^2$.
  Now let $r$ be a prime divisor of $\gcd(q_1,q_2)$ $\Longrightarrow$ $r|h$, $r|i^2+j^2$. With 
  $b=\frac{pi(h^2+j^2)}{q}$ and $\gcd(p,q)=1$ we also have $r|i$ or $r|h^2+j^2$. In the first case we 
  have $r|i^2+j^2$ $\Longrightarrow$ $r|j^2$ $\Longrightarrow$ 
  $r|\gcd(h,i,j)$ $\Longrightarrow$ $r=1$. In the second case we can use $r|h^2+j^2$ and $r|h$ to 
  conclude $r|j^2$. With this and $r|i^2+j^2$ we also get $r|i^2$ and 
  so $r|\gcd(h,i,j)=1$ $\Longrightarrow$ $r=1$. So we know $\gcd(q_1,q_2)$=1.\\[3mm]
  Analog we get $q=q_3q_4$ with $\gcd(q_3,q_4)=1$, $q_3|i$,  and $q_4|\left(h^2+j^2\right)$.\\[3mm]
  Now we set $q_1=w_1w_2$, $q_2=w_3w_4$, $q_3=w_1w_3$, and 
  $q_4=w_2w_4$. With $\gcd(q_1,q_2)=\gcd(q_3,q_4)=1$ we can conclude the 4 divisibility conditions 
  for each $w_i$ and that the $w_i$ are pairwise coprime.
\end{Proof}

\begin{Lemma}
  $$w_4|2(i+h)\,.$$
\end{Lemma}
\begin{Proof}
  We consider $ai-bh=\frac{pih(i+h)(i-h)}{q}$ and conclude $w_4|(i-h)(i+h)$. Now we consider 
  a prime factor $r$ with $r|(i-h)$ and $\gcd(r,i+h)=1$. Because $r|w_4|a,b,c$ we get
  $r|(i^2+j^2)+(h^2+j^2)+2(ih-j^2)=(i+h)^2$, a contradiction to $\gcd(r,i+h)=1$. The proof is completed 
  by $\gcd(i+h,i-h)|2$.
\end{Proof}

\begin{Lemma}
  \label{lem2.3}
  $$w_4\le 8n\,.$$
\end{Lemma}
\begin{Proof}
  To prove the lemma we will show $w_4|8c$. From $w_4|2(i+h)$ we conclude
  $w_4|2(i^2+j^2)+2(h^2+j^2)-2(i+h^2)=4(j^2-ih)$ and thus $w_4|2(i+h)\frac{4(ih-j^2)}{w_4}=8c$.
\end{Proof}

The next step is to find a parameterization of the set of solutions which is better suited for
computational purposes.
Therefore we set $$w_2=st^2$$ and $$w_3=uv^2$$ with squarefree integers $s$ and $u$. Because 
$w_2|h^2+j^2$, $w_2|h$, $w_3|i^2+j^2$, $w_3|i$, and $\gcd(w_2,w_3)=1$ we have $stuv|j$. Thus we can
set
\begin{eqnarray}
  h&=&\alpha w_1st^2,\nonumber\\
 i&=&\beta w_1uv^2,\nonumber\\
  j&=&\gamma stuv\nonumber
\end{eqnarray}
with integers $\alpha$, $\beta$, and $\gamma$.\\
With this we can give the following parameterization of the set of integer Heronian triangles.
\begin{eqnarray}
  a&=&\frac{p\alpha u[(\beta w_1v)^2+(\gamma st)^2]}{w_4},\nonumber\\
 b&=&\frac{p\beta s[(\alpha w_1t)^2+(\gamma uv)^2]}{w_4},\nonumber\\
 c&=&\frac{p(\beta uv^2+\alpha st^2)(\beta\alpha w_1^2-\gamma^2su)}{w_4}.\nonumber
\end{eqnarray}

\section{Algorithms for the generation of integer Heronian triangles}
\label{sec_algo}

\noindent
In this section we list several algorithms to generate all integer Heronian triangles with
diameter at most $n$. The main idea of the first algorithm is to utilize the parameterization of the 
previous section to run through all possible values for $a$, $w_4$ and then to determine all possible
parameters $p$, $w_1$, $s$, $t$, $u$, $v$, $\alpha$, $\beta$ and $\gamma$. Without loss of generality
we can assume that $a\ge b$ and thus $n\le 2a-1$. Then by Lemma \ref{lem2.3} we have $w_4\le 8n\le 16a$.\\

\begin{Algorithm}{\textbf{(Generation of integer Heronian triangles I)}\\}
  \label{algo_i}
  determine the prime factorization of all integers at most $16n$\\
  determine the solutions of $z=x^2+y^2$ for all $z\le 16n$\\
  \texttt{for} $a$ \texttt{from} $1$ \texttt{to} $n$\\
  \hspace*{3mm}\texttt{for} $w_4$ \texttt{from} $1$ \texttt{to} $16a$\\
  \hspace*{6mm}\texttt{loop over} all quadruples $(p,\alpha,u,z)$ with $p\alpha uz$=$aw_4$\\
  \hspace*{9mm}\texttt{loop over} all pairs $(x,y)$ with $x^2+y^2=z$\\
  \hspace*{12mm}\texttt{loop over} all triples $(\beta,w_1,v)$ with $\beta w_1v=x$\\
  \hspace*{15mm}\texttt{loop over} all triples $(\gamma,s,t)$ with $\gamma st=y$\\
  \hspace*{18mm}calculate and output $a$, $b$, $c$\\
\end{Algorithm}

\noindent
In order to prove the running time $\mathcal{O}\left(n^{2+\varepsilon}\right)$ of Algorithm \ref{algo_i} we rephrase two results from number theory.

\begin{Theorem}{(\textbf{Theorem 317 \cite{number_theory_hw}})}
  \label{lemma_num_theo}
  For $\varepsilon>0$ and $n>n_0(\varepsilon)$
  $$\tau(n)<2^{(1+\varepsilon)\frac{\log n}{\log\log n}}$$
  where $\tau(n)$ denotes the number of divisors of $n$.
\end{Theorem}

\noindent
So for each $\varepsilon>0$, $f\le 16n^2$ there are only $\mathcal{O}\left(n^\varepsilon\right)$ quadruples
$(f_1,f_2,f_3,f_4)$ with $f_1f_2f_3f_4=f$.

\begin{Lemma}
  The equation $z=x^2+y^2$ has at most $\mathcal{O}\left(z^\varepsilon\right)$ solutions in positive integers $x$, $y$ for
  each $\varepsilon>0$.
\end{Lemma}
\begin{Proof}
  If we denote the number of solutions of $z=x^2+y^2$ in pairs $(x,y)$ of integers by $r_2(z)$ then we
  have \cite{0086.26202,sum_of_squares}
  $$
    r_2(z)=4\cdot\sum_{d|z}\sin\left(\frac{1}{2}\pi d\right)\in \mathcal{O}\left(z^\varepsilon\right).
  $$
\end{Proof}

\noindent
Thus for each $\varepsilon>0$ and each $z\le 16n^2$ there are only $\mathcal{O}\left(n^\varepsilon\right)$ integer 
solutions of $z=x^2+y^2$. Consequently there exists an implementation of Algorithm \ref{algo_i} with running time
$\mathcal{O}\left(n^{2+\varepsilon}\right)$. Furthermore we can conclude that there are $\mathcal{O}\left(n^{1+\varepsilon}\right)$
integer Heronian triangles with diameter $n$. Maybe a faster algorithm can be designed by using refined number 
theoretic conditions on $w_4$. Unfortunately we were not able to find estimations on the number of integer
Heronian triangles in literature. Therefore we are unable to give a lower bound for the complexity of generating
integer Heronian triangles.

In order to derive a second algorithm for the determination of integer Heronian triangles we utilize the Heron formula 
for the area of a triangle $\Delta=(a,b,c)$ and consider 
$$
  16A^2=(p-c)(p+c)(c-q)(c+q)
$$ 
with $p=a+b$ and $q=a-b$.

The idea is to run trough all possible values for $4A$ and then determine $a$, $b$ and $c$ by factorising $16A^2$.

\begin{Algorithm}
 {\textbf{(Generation of integer Heronian triangles II)}\\} 
 \label{algo_ii}
 \texttt{loop over} all $m$ and the prime factorization of $m^2$ with $1\le m\le\sqrt{3}n^2$\\ 
 \hspace*{3mm}\texttt{loop over} all $p-c$, $p+c$, $c-q$, $c+q$ with $m^2=(p-c)(p+c)(c-q)(c+q)$\\
 \hspace*{6mm}determine $a$, $b$, and $c$\\
 \hspace*{9mm}\texttt{if} $a$, $b$, and $c$ are positive integers satisfying the triangle conditions
              \texttt{then} output $a$, $b$, and $c$
\end{Algorithm}

\noindent
Since $16A^2=(a+b+c)(a+b-c)(a-b+c)(-a+b+c)\le 3n^4$ we have $m=4A\le\sqrt{3}n^2$. For the factorization 
of $m$ we may use an arbitrary algorithm with running time $\mathcal{O}\left(m^\varepsilon\right)$ \cite{0508.10004}.
If we are allowed to use $\Omega(n^2)$ space a less sophisticated possibility would be to use the Sieve of 
Eratosthenes on the numbers $1$ to $\sqrt{3}n^2$. Thus Algorithm \ref{algo_ii} can be implemented with running time $\mathcal{O}\left(n^{2+\varepsilon}\right)$.

For completeness we would also like to give the pseudo code of the algorithm mentioned in the introduction.

\begin{Algorithm}
 {\textbf{(Generation of integer Heronian triangles III)}\\}
 \label{algo_iii}
  \texttt{for} $a$ \texttt{from} $1$ \texttt{to} $n$\\
 \hspace*{3mm}\texttt{for} $b$ \texttt{from} $\left\lceil\frac{a+1}{2}\right\rceil$ \texttt{to} $a$\\
 \hspace*{6mm}\texttt{for} $c$ \texttt{from} $a+1-b$ \texttt{to} $b$\\
 \hspace*{9mm}\texttt{if} $(a+b+c)(a+b-c)(a-b+c)(-a+b+c)$ is the square of an integer\\
 \hspace*{9mm}\texttt{then} output $a$, $b$, and $c$
\end{Algorithm}

\noindent
The running time of Algorithm \ref{algo_iii} is $\mathcal{O}\left(n^3\right)$. It has the advantage of producing only one representative of each equivalence class of integer Heronian triangles in a canonical ordering. Due to the overhead of Algorithm \ref{algo_i} and Algorithm \ref{algo_ii} the trivial Algorithm \ref{algo_iii} is faster for \textit{small} values of $n$.

For a practical implementation we describe some useful tricks to enhance Algorithm \ref{algo_iii} a bit. We observe that if $(a+b+c)(a+b-c)(a-b+c)(-a+b+c)$ is a square then it must also be a square if we calculate in the ring $\mathbb{Z}_m$ for all $m\in\mathbb{N}$. In our implementation we have used the set of divisors of $420$ for $m$. In a precalculation we have determined
all possible triples in $\mathbb{Z}_{420}^3$. Hereby the number of candidates is reduced by a factor of $\frac{14744724}{74088000}\approx 0.199$. Additionally we determine the squarefree parts of the integers at most $3n$ in a precalculation. Instead of determining the square root of a big integer we determine squarefree parts of integers. If $\text{sfp}(f)$ gives the squarefree part of $f$ then we have $$\text{sfp}(f_1\cdot f_2)=\frac{\text{sfp}(f_1)\cdot\text{sfp}(f_2)}{\gcd\left(\text{sfp}(f_1),\text{sfp}(f_2)\right)^2}.$$
Thus we can avoid high precision arithmetic by using a $\gcd$-algorithm. Without it we would have to deal with very large numbers -- since we compute up to $n=600000$. A complete list of the integer Heronian triangles of diameter at most $600000$ can be obtained upon request to the author. In the following sections we will use this list to attack several combinatorial problems.

%\section{Conclusion}
%Now we would like to compare the practical running times of the described algorithms.

%\begin{table}[h]
%\caption[Comparison of running times]{Comparison of running times}
% \begin{tabular}{r|r|r|r}
%  $n$ & Algorithm \ref{algo_i} & Algorithm \ref{algo_ii} & Algorithm \ref{algo_iii}  \\ 
% \hline
%  100 & 2.67 s  & 1.08 s & 0.02 s \\ 
%  200 & 17.0 s  & 4.85 s & 0.11 s \\ 
%  400 & 104 s   & 25.7 s & 0.80 s \\ 
%  800 & 620 s   & 140 s  & 6.47 s \\ 
% 1600 & 3599 s  & 756 s  & 51.2 s \\ 
% 3200 & 20403 s & 3987 s & 403 s  \\ 
%\end{tabular} 
%\end{table}
% 
% If we extrapolate the running times then we can state that for $n\ge 80,000$ Algorithm
% \ref{algo_i} beats Algorithm \ref{algo_iii}. 

\section{Maximal integral triangles}

\noindent
A result due to Almering \cite{rationale_vierecke} is the following. Given any rational triangle $\Delta=(a,b,c)\in\mathbb{Q}^3$ in the plane, i.~e. a triangle with rational side lengths, the set of all points $x$ with rational distances to the three corners of $\Delta$ is dense in the plane. Later Berry \cite{berry} generalized this results to triangles which side lengths are rational when squared and with one side length rational. If we proceed to integral side lengths and integral coordinates the situation is a bit different. In \cite{paper_carpet} the authors search for inclusion-maximal integral triangles over $\mathbb{Z}^2$ and answer the existence question from \cite{gauss_integers} positively. They exist but appear to be somewhat rare. There are only seven inclusion-maximal integral triangles with diameter at most $5,000$.

{\sloppy
Here we have used the same algorithm as in \cite{paper_carpet} to determine inclusion-maximal integral triangles over $\mathbb{Z}^2$ with diameter at most $15,000$. Up to symmetry the complete list is given by:\\
\small{
(2066,1803,505),  (2549,2307,1492), (3796,2787,2165), (4083,2425,1706), (4426,2807,1745), (4801,2593,2210), \\
(4920,4177,985),  (5044,4443,2045), (5045,4803,244),  (5186,5163,745),  (5905,5763,1586), (5956,4685,2427), \\
(6120,5953,409),  (6252,3725,3253), (6553,5954,3099), (6577,5091,1586), (6630,5077,1621), (6787,5417,1546), \\
(6855,6731,130),  (6890,6001,1033), (6970,4689,4217), (6987,5834,1585), (7481,6833,5850), (7574,4381,3207), \\
(7717,6375,1396), (7732,7215,541),  (7734,6895,4537), (7793,4428,3385), (7837,6725,1308), (7913,6184,1745), \\
(7985,7689,298),  (8045,7131,1252), (8187,6989,1252), (8237,7899,4036), (8249,7772,879),  (8286,5189,3865), \\
(8375,6438,1949), (8425,4706,3723), (8644,7995,1033), (8961,8633,740),  (9683,8749,4632), (9745,5043,4706), \\
(9771,7373,5044),   (9840,8473,2089),   (9939,6388,3845),   (9953,6108,4825),   (10069,9048,6421),   \\
(10081,8705,1378),  (10088,8886,4090),  (10090,9606,488),   (10100,5397,5389),  (10114,5731,4405),   \\
(10372,7739,2775),  (10394,8499,1993),  (10441,6122,5763),  (10595,10283,340),  (10600,6737,3881),   \\
(10605,8957,1754),  (10615,10119,562),  (10708,9855,1069),  (10804,8691,7013),  (10825,8259,3242),   \\
(10875,9805,1076),  (10993,8164,3315),  (11133,10250,6173), (11199,10444,757),  (11283,8788,4229),   \\
(11332,9147,6029),  (11434,6159,5305),  (11441,7577,3880),  (11559,6145,5416),  (11765,10892,877),   \\
(11787,9341,3172),  (12053,8979,3076),  (12676,9987,3845),  (12745,12603,1586), (12757,11544,1237),  \\
(12810,12077,2669), (12818,11681,1601), (12946,9523,3425),  (12953,8361,4930),  (12965,12605,5406),  \\
(13012,11405,2091), (13061,9745,8934),  (13100,12875,1011), (13106,11908,6198), (13115,11492,1709),  \\
(13130,12097,2329), (13309,12916,8585), (13350,7901,5645),  (13369,12867,698),  (13385,11931,1618),  \\
(13445,9750,3701),  (13466,8665,4803),  (13683,8042,6841),  (13700,11115,2621), (13710,13462,260),   \\
(13740,8053,5951),  (13780,12002,2066), (13876,10657,3315), (13940,9378,8434),  (13940,13647,12775), \\
(13951,11785,9608), (13971,10804,8933), (14065,10984,3831), (14065,12531,5378), (14126,12135,4357),  \\
(14172,12725,1933), (14185,9879,8194),  (14282,8665,5619),  (14331,7517,6964),  (14356,14019,2837),  \\
(14379,10685,3748), (14545,13274,2427), (14615,11332,3441), (14625,13060,1799), (14633,12329,2320),  \\
(14637,11170,6641), (14677,11436,6145), (14690,11353,5641), (14700,11861,2845), (14705,14351,8482),  \\
(14775,10673,5284), (14785,12219,3242), (14819,13810,1695), and (14962,13666,11700).
}}

Thus with $126$ examples the situation changes a bit. There do exist lots of inclusion-maximal integral triangles over $\mathbb{Z}^2$. Some triangles of this list may be derived from others since $\left(\frac{a}{g},\frac{b}{g},\frac{c}{g}\right)$ is an inclusion-maximal integral triangle over $\mathbb{Z}^2$ for $g=\gcd(a,b,c)$ if $(a,b,c)$ is an inclusion-maximal integral triangle over $\mathbb{Z}^2$. Here the limiting factor is the algorithm from \cite{paper_carpet} and not the generation of integral Heronian triangles. We remark that there are also inclusion-maximal integral tetrahedrons over $\mathbb{Z}^3$ \cite{paper_carpet}.

\section{$\mathbf{n_2}$-cluster}

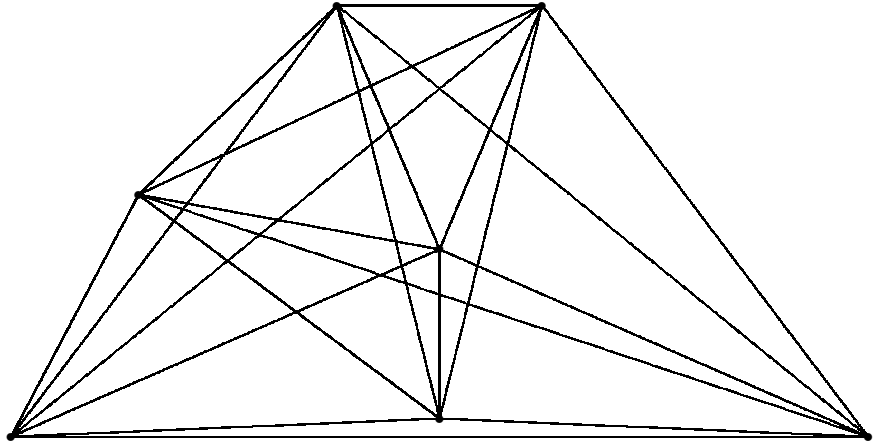
\begin{figure}[h]
  \begin{center}
    \setlength{\unitlength}{0.0003mm}
    \begin{picture}(376000,191000)
      \put(0,0){\circle*{3000}}
      \put(375360,0){\circle*{3000}}
      \put(55860,106855){\circle*{3000}}
      \put(187680,7990){\circle*{3000}}
      \put(187680,82688){\circle*{3000}}
      \put(142800,190400){\circle*{3000}}
      \put(232560,190400){\circle*{3000}}
      \qbezier(0,0)(187680,0)(375360,0)
      \qbezier(0,0)(27930,53427.5)(55860,106855)
      \qbezier(0,0)(93840,3995)(187680,7990)
      \qbezier(0,0)(93840,41344)(187680,82688)
      \qbezier(0,0)(71400,95200)(142800,190400)
      \qbezier(0,0)(116280,95200)(232560,190400)
      \qbezier(375360,0)(215610,53427.5)(55860,106855)
      \qbezier(375360,0)(281520,3995)(187680,7990)
      \qbezier(375360,0)(281520,41344)(187680,82688)
      \qbezier(375360,0)(259080,95200)(142800,190400)
      \qbezier(375360,0)(303960,95200)(232560,190400)
      \qbezier(55860,106855)(121770,57422.5)(187680,7990)
      \qbezier(55860,106855)(121770,94771.5)(187680,82688)
      \qbezier(55860,106855)(99330,148627.5)(142800,190400)
      \qbezier(55860,106855)(144210,148627.5)(232560,190400)
      \qbezier(187680,7990)(187680,45339)(187680,82688)
      \qbezier(187680,7990)(165240,99195)(142800,190400)
      \qbezier(187680,7990)(210120,99195)(232560,190400)
      \qbezier(187680,82688)(165240,136544)(142800,190400)
      \qbezier(187680,82688)(210120,136544)(232560,190400)
      \qbezier(142800,190400)(187680,190400)(232560,190400)
    \end{picture}
  \end{center}
  \caption{Almost a $7_2$-cluster.}
  \label{fig_almost}
\end{figure}

\noindent
A $n_2$ cluster is a set of $n$ lattice points in $\mathbb{Z}^2$ where all pairwise distances are integral, no three points are on a line, and no four points are situated on a circle \cite{cluster}. The existence of a $7_2$-cluster is an unsolved problem of \cite[Problem D20]{upin} and \cite{cluster}. Since a $7_2$-cluster is composed of Heronian triangles and a special case of a plane integral point set, we can use the exhaustive generation algorithms described in \cite{phd_kurz,paper_alfred} to search for $7_2$-clusters. The point set with coordinates
\begin{eqnarray*}
  &&\Big\{(0,0), (375360,0), (55860,106855), (187680,7990),\\&& (187680,82688),(142800,190400), (232560,190400)\Big\}
\end{eqnarray*}
and distance matrix
$$
  \left(
    \begin{array}{rrrrrrr}
           0 & 375360 & 120575 & 187850 & 205088 & 238000 & 300560 \\ 
      375360 &      0 & 336895 & 187850 & 205088 & 300560 & 238000 \\ 
      120575 & 336895 &      0 & 164775 & 134017 & 120575 & 195455 \\ 
      187850 & 187850 & 164775 &      0 &  74698 & 187850 & 187850 \\ 
         205088 & 205088 & 134017 &  74698 &      0 & 116688 & 116688 \\
         238000 & 300560 & 120575 & 187850 & 116688 &      0 &  89760 \\
      300560 & 238000 & 195455 & 187850 & 116688 &  89760 &      0
    \end{array}
  \right)
$$
is an integral point set over $\mathbb{Z}^2$, see Figure \ref{fig_almost}. Unfortunately the points $1$, $2$, $6$ and $7$ are on a circle. But, no three points are on a line and no other quadruple is
situated on a circle. If we add $(319500,106855)$ as an eighth point we receive an integral point set $\mathcal{P}$ 
over $\mathbb{Z}^2$ where no three points are situated on a line. There are exactly three quadruples of points which are 
situated on a circle: $\{1,2,3,8\}$, $\{1,2,6,7\}$, and $\{3,6,7,8\}$.

\noindent
We would like to remark that the automorphism group of an $n_2$-cluster for $n\ge 6$ must be trivial. In \cite{phd_kurz} the possible automorphism groups of planar integral point sets were determined to be isomorphic to $\text{id}$, $\mathbb{Z}_2$, $\mathbb{Z}_2\times\mathbb{Z}_2$, $\mathbb{Z}_3$ or $S_3$. It was also shown that an automorphism of order $3$ is only possible for characteristic $3$. Since $n_2$-cluster have characteristic $1$ such an automorphism cannot exist. If we would have an automorphism of order $2$ then for $n\ge 6$ either three points are collinear or four points are situated on a circle.

Using our list of Heronian triangles we have performed an exhaustive search for $7_2$-clusters up to diameter $600000$, unfortunately without success. If we relax the condition of integral coordinates, then examples do exist, see \cite{paper_kreisel}.

\section{Perfect pyramids}

\begin{figure}[htp]
  \begin{center}
    \setlength{\unitlength}{3.cm}
    \begin{picture}(1,0.9)
      \put(0,0){\line(1,0){1}}
      \qbezier(0,0)(0.25,0.433)(0.5,0.866)
      \qbezier(1,0)(0.75,0.433)(0.5,0.866)
      \put(0.5,0.866){\line(0,-1){0.566}}
      \qbezier(0,0)(0.25,0.15)(0.5,0.3)
      \qbezier(1,0)(0.75,0.15)(0.5,0.3)
      \put(0.15,0.4){$c$}
      \put(0.80,0.4){$b$}
      \put(0.48,0.02){$a$}
      \put(0.52,0.45){$d$}
      \put(0.27,0.2){$e$}
      \put(0.70,0.22){$f$}
    \end{picture}
    \caption{The six edges of a tetrahedron $(a,b,c,d,e,f)$.}
    \label{fig_tetrahedron}
  \end{center}
\end{figure}
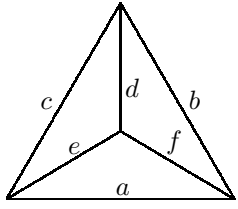

\noindent
In \cite{perfect_pyramids} the author considers tetrahedra with integral side lengths, integral face areas, and integral volume, see also \cite{petersen}. The smallest such example has side lengths $(a,b,c,d,e,f)=(117,84,51,52,53,80)$ using the notation from Figure \ref{fig_tetrahedron}. In the plane a triangle with integral edge lengths and rational area is forced to have an integral area. The situation changes slightly in three-dimensional space. Here it is possible that the edge lengths of a tetrahedron are integral and that the volume is genuinely rational. If the edge lengths are integral and the face areas and the volume are rational, then all values must be integral, see \cite{jnt_2006}. In \cite{perfect_pyramids} it was also shown that a perfect pyramid with at most two different edge lengths cannot exist. For three different edge lengths a parameter solution of an infinite family is given.

\noindent
More comprising the authors of \cite{perfect_pyramids,jnt_2006} consider all possibilities of coincidences of edge lengths. Up to symmetry we have the following configurations:
$$
\begin{array}{lrlrl}
  \text{1-parameter} &\text{(i)}   &\!\!\!\!a=b=c=d=e=f,                                      \\
  \text{2-parameter} &\text{(i)}   &\!\!\!\!a=b=c=d=e,f, &\text{(ii)}   &\!\!\!\! a=b=c=d,e=f,\\
                     &\text{(iii)} &\!\!\!\!a=c=d=f,b,e, &\text{(iv)}   &\!\!\!\! a=b=c,d=e=f,\\
                     &\text{(v)}   &\!\!\!\!a=d=f,b=c=e,                                      \\
  \text{3-parameter} &\text{(i)}   &\!\!\!\!a=b=c=d,e,f, &\text{(ii)}   &\!\!\!\! a=c=d=f,b,e,\\
                     &\text{(iii)} &\!\!\!\!a=b=c,d=e,f, &\text{(iv)}   &\!\!\!\! a=d=f,b=c,e,\\
                     &\text{(v)}   &\!\!\!\!a=d=f,b=e,c, &\text{(vi)}   &\!\!\!\! a=d,b=e,c=f,\\
                     &\text{(vii)} &\!\!\!\!a=e,b=f,c=d, &\text{(viii)} &\!\!\!\! a=b,c,d=e=f,\\
                     &\text{(ix)}  &\!\!\!\!a=d,b=f,c=e,                                      \\
  \text{4-parameter} &\text{(i)}   &\!\!\!\!a=b=c,d,e,f, &\text{(ii)}   &\!\!\!\! a=b=d,c,e,f,\\
                     &\text{(iii)} &\!\!\!\!a=b=f,c,d,e, &\text{(iv)}   &\!\!\!\! a=d,b=e,c,f,\\
                     &\text{(v)}   &\!\!\!\!a=d,b=f,c,e, &\text{(vi)}   &\!\!\!\! a=b,d=f,c,e,\\
                     &\text{(vii)} &\!\!\!\!a=b,d=e,c,f,                                      \\
  \text{5-parameter} &\text{(i)}   &\!\!\!\!a=d,b,c,e,f, &\text{(ii)}   &\!\!\!\! a=b,c,d,e,f,\\
  \text{6-parameter} &\text{(i)}   &\!\!\!\!a,b,c,d,e,f.                                      \\
\end{array}
$$
In the cases 1(i), 2(i-v), 3(i), 3(iii), 3(iv), 3(vii), and 4(i) there is no solution. There are parameter solutions of infinite families known for the cases 3(ii), 3(v), 3(vi), 3(viii), and 3(ix). Unfortunately in case 3(viii) all known solutions are degenerated, meaning that the corresponding pyramid has volume $0$. It is a conjecture of \cite{jnt_2006} that no non-degenerated solution exists in this case. For the cases 4(iv), 4(vii), 5(ii), and 6(i) sporadic solutions are known. Cases 4(ii), 4(iii), 4(v), 4(vi), and 5(i) remain open problems.

In order to answer some of these questions we have performed an exhaustive search on perfect pyramids up to diameter $600000$ using the following algorithm:
Let $n$ be the maximum diameter, $\kappa\in \mathbb{N}$, and $\varphi:\{1,\dots,n\}\times\{1,\dots,n\}\rightarrow\{0,\kappa-1\}$ a mapping. For $0\le i\le \kappa-1$ let $\mathcal{L}_i$ contain integers $1\le c\le n$ such that there exit integers $a,b\le n$ fulfilling $\phi(a,b)=i$ and where $(a,b,c)$ is a Heronian triangle. 
\begin{Algorithm}{\textbf{(Generation of perfect pyramids)}\\}
  \label{algo_pp}
  \texttt{for} $d$ \texttt{from} $1$ \texttt{to} $n$\\
  \hspace*{3mm}\texttt{loop over} all $\left(a_1,b_1\right)$, $\left(a_2,b_2\right)$ where $\left(d,a_1,b_1\right)$
               and $\left(d,a_2,b_2\right)$ are Heronian triangles\\
  \hspace*{6mm}$lb=\max\Big\{\left|a_1-a_2\right|,\left|b_1-b_2\right|\Big\}+1$\\
  \hspace*{6mm}$ub=\min\left\{a_1+a_2-1,b_1+b_2-1,d\right\}$\\
  \hspace*{6mm}\texttt{if} $ub-lb+1=\min\Big\{ub-lb+1,\left|\mathcal{L}_{\varphi(a_1,a_2)}\right|,
               \left|\mathcal{L}_{\varphi(b_1,b_2)}\right|\Big\}$ \texttt{then}\\
  \hspace*{9mm}\texttt{for} $x$ \texttt{from} $lb$ \texttt{to} $ub$\\
  \hspace*{12mm}\texttt{if} $P=(d,a_1,a_2,b_1,b_2,x)$ is a perfect pyramid \texttt{then} output $P$\\
  \hspace*{6mm}\texttt{if} $\left|\mathcal{L}_{\varphi(a_1,a_2)}\right|<\min\Big\{ub-lb+1,
               \left|\mathcal{L}_{\varphi(b_1,b_2)}\right|\Big\}$ \texttt{then}\\
  \hspace*{9mm}\texttt{for} $j$ \texttt{from} $1$ \texttt{to} $\left|\mathcal{L}_{\varphi(a_1,a_2)}\right|$\\
  \hspace*{12mm}$x=\mathcal{L}_{\varphi(a_1,a_2)}(j)$\\
  \hspace*{12mm}\texttt{if} $P=(d,a_1,a_2,b_1,b_2,x)$ is a perfect pyramid \texttt{then} output $P$\\
  \hspace*{6mm}\texttt{if} $\left|\mathcal{L}_{\varphi(b_1,b_2)}\right|<\min\Big\{ub-lb+1,
               \left|\mathcal{L}_{\varphi(a_1,a_2)}\right|\Big\}$ \texttt{then}\\
  \hspace*{9mm}\texttt{for} $j$ \texttt{from} $1$ \texttt{to} $\left|\mathcal{L}_{\varphi(b_1,b_2)}\right|$\\
  \hspace*{12mm}$x=\mathcal{L}_{\varphi(b_1,b_2)}(j)$\\
  \hspace*{12mm}\texttt{if} $P=(d,a_1,a_2,b_1,b_2,x)$ is a perfect pyramid \texttt{then} output $P$\\
\end{Algorithm}

\noindent
The efficiency of Algorithm \ref{algo_pp} depends on a suitable choice of $\kappa$ and $\varphi$ in order to keep the lists $\mathcal{L}_i$ small. From a theoretical point of view for given integers $a$ and $b$ there do exist at most $4\cdot\tau(ab)^2$ different values $c$ such that $(a,b,c)$ is a Heronian triangle \cite{two_fixed}. Here $\tau(m)$ denotes the number of divisors of $m$ and we have $\tau(m)\in\mathcal{O}\left(m^\varepsilon\right)$ for all $\varepsilon>0$.

Unfortunately we have found no examples for one of the open cases. Thus possible examples have a diameter greater than $600000$.

In \cite{jnt_2005} the authors have considered rational tetrahedra with edges in arithmetic progression. They proved that tetrahedron with integral edge lengths, rational face areas and rational volume do not exist. If only one face area is forced to be rational then there exist the example $(a,b,c,d,e,f)=(10,8,6,7,11,9)$, which is conjectured to be unique up to scaling. We have verified this conjecture up to diameter $600000$.

\begin{table}[htp]
  \begin{center}
    \begin{tabular}{|r|r|r|r|r|r|r|r|}
      \hline
      \textbf{surface area} & \textbf{volume} & \textbf{a} & \textbf{b} & \textbf{c} & \textbf{d} & \textbf{e} & \textbf{f} \\
      \hline
      6384 & 8064 & 160 & 153 & 25 & 39 & 56 & 120 \\
      \hline
      64584 & 170016 & 595 & 208 & 429 & 116 & 325 & 276 \\
      64584 & 200928 & 595 & 116 & 507 & 208 & 325 & 276 \\
      \hline
      69058080 & 14985432000 & 11660 & 5512 & 10812 & 1887 & 9945 & 5215 \\
      69058080 & 15020132400 & 11687 & 5215 & 10812 & 1887 & 9945 & 5512 \\
      69058080 & 16198182000 & 11687 & 1887 & 11660 & 5215 & 9945 & 5512 \\
      \hline
    \end{tabular}
    \caption{Sets of primitive perfect pyramids with equal surface area.}
    \label{table_equal_surface_area}
  \end{center}
\end{table}

\noindent
Now we consider sets of primitive perfect pyramids (here the greatest common divisor of the six edge length must be equal to one) which have equal surface area, see \cite{petersen}. We have performed an exhaustive search on the perfect pyramids up to diameter $600000$ and list the minimal sets, with respect to the surface area, in Table \ref{table_equal_surface_area}. We would like to remark that six triples with equal surface area were found. Clearly the same question arises for equal volume instead of equal surface area. In \cite{petersen} the smallest examples for up to three pyramids with equal volume are given. Unfortunately we have found no further examples consisting of three or more primitive perfect pyramids up to diameter $600000$. Additionally we have searched, without success, for a pair of perfect pyramids with equal surface area and equal volume.

% 3 -> 1037897952 518266148400
% 3 -> 1054812360 1000867140000
% 3 -> 122084820 25059434400
% 3 -> 65396173920 602133130368000
% 3 -> 69058080 14985432000
% 3 -> 77934748320 576889369056000

Clearly the concept of perfect pyramids can be generalized to higher dimensions. The $m$-dimensional volume $V_m\left(\mathcal{S}\right)$ of a simplex $\mathcal{S}=(v_0,\dots,v_m)$ can be expressed by a determinant \cite{volumen}. Therefore let $D:=\Big(\left\Vert v_i-v_j\right\Vert_2^2\Big)_{0\le i,j\le m}$. With this we have $V_m\left(\mathcal{S}\right)^2=\frac{(-1)^{m+1}}{2^m(m!)^2}\cdot\det B$, where $B$ arises from $D$ by bordering $D$ with a top row $(0,1,\dots,1)$ and a left column $(0,1,\dots,1)^T$. We call a $m$-dimensional simplex $\mathcal{S}$ with integral edge lengths a perfect pyramid if $V_r\Big(\left\{i_0,\dots,i_r\right\}\Big)\in\mathbb{Q}_{>0}$ for all $2\le r\le m$ and all $\{i_0,\dots,i_r\}\subseteq\mathcal{S}$. Up to diameter $600000$ we have only found some degenerated perfect pyramids where $V_4\left(\mathcal{S}\right)=0$.

\section{Heronian triangles and sets of Heronian triangles with special properties}

\noindent
In \cite{van_luijk} it is shown that for every positive integer $N$ there exists an infinite family parameterized by $s\in\mathbb{Z}_{>0}$, of $N$-tuples of pairwise nonsimilar Heron triangles, all $N$ with the same area $A(s)$ and the same perimeter $p(s)$, such that for any two different $s$ and $s'$ the corresponding ratios $A(s)/p(s)^2$ and $A(s')/p(s')^2$ are different. Randall Rathbun found the smallest $N$-tuples for $N\le 9$. In Table \ref{table_tuples_perimeter_area} we give the perimeter and the area of the smallest $N$-tuples for $N\le 10$. The smallest $11$-tuple has a perimeter greater than $1200000$.

\begin{table}[htp]
  \begin{center}
    \begin{tabular}{|r|r|r||r|r|r|}
      \hline
      $N$ & perimeter & area & $N$ & perimeter & area \\
      \hline
        1 &     12 &          24 &  6 &   2340 &      786240 \\
        2 &     70 &         840 &  7 &  11700 &    19656000 \\
        3 &     98 &        1680 &  8 &  84630 &   365601600 \\
        4 &    448 &       26880 &  9 & 142912 &  2117955840 \\
        5 &   1170 &      196560 & 10 & 441784 & 30324053760 \\
      \hline
    \end{tabular}\\[2mm]
    \caption{Minimum perimeter of $N$-tuples of Heronian triangles with equal perimeter and equal area.}
    \label{table_tuples_perimeter_area}
  \end{center}
\end{table}

\noindent
Some authors consider Heronian triangles with rational medians, see e.~g. \cite{0873.11022,0923.11052}. Infinite families and some sporadic examples of Heronian triangles with two rational medians are know. Whether there exists a Heronian triangle with three rational medians is an open question. For this problem we can only state that the lists in \cite{0873.11022,0923.11052} are complete up to diameter $600000$.

\nocite{perfect_pyramids,jnt_2006,jnt_2005,upin,dove}

\bibliography{On_Heronian_Triangles}

\begin{thebibliography}{10}

\bibitem{rationale_vierecke}
J.H.J. Almering.
\newblock {Rational quadrilaterals}.
\newblock {\em Nederl. Akad. Wet., Proc., Ser. A}, 66:192--199, 1963.

\bibitem{berry}
T.G. Berry.
\newblock Points at rational distance from the vertices of a triangle.
\newblock {\em Acta Arith.}, 62(4):391--398, 1992.

\bibitem{perfect_pyramids}
R.H. Buchholz.
\newblock Perfect pyramids.
\newblock {\em Bull. Aust. Math. Soc.}, 45(3):353--368, 1992.

\bibitem{0873.11022}
R.H. Buchholz and R.L. Rathbun.
\newblock An infinite set of heron triangles with two rational medians.
\newblock {\em Am. Math. Mon.}, 104(2):107--115, 1997.

\bibitem{0923.11052}
R.H. Buchholz and R.L. Rathbun.
\newblock Heron triangles and elliptic curves.
\newblock {\em Bull. Aust. Math. Soc.}, 58(3):411--421, 1998.

\bibitem{0219.10027}
J.R. Carlson.
\newblock Determination of heronian triangles.
\newblock {\em Fibonacci Q.}, 8:499--506, 1970.

\bibitem{jnt_2005}
C.~Chisholm and J.A. MacDougall.
\newblock Rational tetrahedra with edges in arithmetic progression.
\newblock {\em J. Number Theory}, 111(1):57--80, 2005.

\bibitem{jnt_2006}
C.~Chisholm and J.A. MacDougall.
\newblock Rational and heron tetrahedra.
\newblock {\em J. Number Theory}, 121(1):153--185, 2006.

\bibitem{gauss_integers}
S.~Dimiev and K.~Markov.
\newblock Gauss {I}ntegers and {D}iophantine {F}igures.
\newblock {\em Mathematics and Mathematical Education}, 31:88--95, 2002.
\newblock arXiv:math.NT/0203061v1 7 Mar 2002.

\bibitem{dove}
K.L. Dove and J.S. Sumner.
\newblock Tetrahedra with integer edges and integer volume.
\newblock {\em Math. Mag.}, 65(2):104--111, 1992.

\bibitem{upin}
R.K. Guy.
\newblock {\em Unsolved problems in number theory. 3rd ed.}
\newblock Problem Books in Mathematics. New York, NY: Springer-Verlag. xviii,
  437~p., 2004.

\bibitem{0086.26202}
G.H. Hardy.
\newblock {\em Ramanujan. Twelve lectures on subjects suggested by his life and
  work}.
\newblock New York: Chelsea Publishing Company. 236 p., 1959.

\bibitem{number_theory_hw}
G.H. Hardy and E.M. Wright.
\newblock {\em An introduction to the theory of numbers}.
\newblock Oxford etc.: Oxford at the Clarendon Press. XVI, 426 p., 5th edition,
  1979.

\bibitem{two_fixed}
E.J. Ionascu, F.~Luca, and P.~St{\u a}nic{\u a}.
\newblock Heron triangles with two fixed sides.
\newblock {\em J. Number Theory}, 126(1):52--67, 2007.

\bibitem{paper_carpet}
A.~Kohnert and S.~Kurz.
\newblock A note on {E}rd\"os-{D}iophantine graphs and {D}iophantine carpets.
\newblock {\em Mathematica Balkanica}, 21(1-2):1--5, 2007.

\bibitem{paper_kreisel}
T.~Kreisel and S.~Kurz.
\newblock There are integral heptagons, no three points on a line, no four on a
  circle.
\newblock {\em Discrete. Comput. Geom.}, 39(4):786--790, 2008.

\bibitem{phd_kurz}
S.~Kurz.
\newblock {\em Konstruktion und {E}igenschaften ganzzahliger {P}unktmengen}.
\newblock PhD thesis, Bayreuth. Math. Schr. 76. Universit\"at Bayreuth, 2006.

\bibitem{paper_alfred}
S.~Kurz and A.~Wassermann.
\newblock On the minimum diameter of plane integral point sets.
\newblock {\em Ars Combin.}, 101:265--287, 2011.

\bibitem{cluster}
L.C. Noll and D.I. Bell.
\newblock $n$-clusters for $1\le n\le 7$.
\newblock {\em Math. Comput.}, 53(187):439--444, 1989.

\bibitem{petersen}
I.~Peterson.
\newblock Perfect pyramids.
\newblock http://www.sciencenews.org/articles/20030726/mathtrek.asp, 2003.

\bibitem{0508.10004}
C.~Pomerance.
\newblock Analysis and comparison of some integer factoring algorithms.
\newblock In H.W. Lenstra and R~Tijdeman, editors, {\em Computational methods
  in number theory}, number 154 in Mathematical {C}entre tracts, pages 89--139.
  Mathematisch {C}entrum, 1982.

\bibitem{0893.11009}
D.J. Rusin.
\newblock Rational triangles with equal area.
\newblock {\em New York J. Math.}, 4:1--15, 1998.

\bibitem{0986.51030}
K.R.S. Sastry.
\newblock Heron triangles: A gergonne-cevian-and-median perspective.
\newblock {\em Forum Geom.}, 1:17--24, 2001.

\bibitem{0255.10019}
D.~Singmaster.
\newblock Some corrections to {C}arlson's "determination of heronian
  triangles".
\newblock {\em Fibonacci Q.}, 11:157--158, 1973.

\bibitem{volumen}
D.M.Y. Sommerville.
\newblock {\em An introduction to the geometry of $n$ dimensions.}
\newblock New York: Dover Publications, Inc. XVII, 196 p., 1958.

\bibitem{van_luijk}
R.~Van~Luijk.
\newblock An elliptic {K}3 surface associated to {H}eron triangles.
\newblock {\em J. Number Theory}, 123(1):92--119, 2007.

\bibitem{sum_of_squares}
E.W. Weisstein.
\newblock Sum of squares function.
\newblock http://mathworld.wolfram.com/SumofSquaresFunction.html, 2003.

\bibitem{0937.11012}
P.~Yiu.
\newblock Construction of indecomposable heronian triangles.
\newblock {\em Rocky Mt. J. Math.}, 28(3):1189--1202, 1998.

\end{thebibliography}
\bibdata{On_Heronian_Triangles}
\bibliographystyle{plain} 

\end{document}